\documentclass[a4paper, 10pt, twosides]{amsart}
\usepackage{lmodern}
\usepackage{amsmath,amsthm}
\usepackage{amssymb}

\usepackage{mathtools}
\usepackage{bbm}
\usepackage{enumitem}
\setenumerate{nolistsep}
\setitemize{nolistsep}
\usepackage[norelsize, ruled, vlined]{algorithm2e}
\usepackage[usenames, dvipsnames]{color}

\usepackage{tikz}
\usetikzlibrary{decorations.markings}
\usetikzlibrary{shapes.geometric,positioning}
\usepackage[utf8]{inputenc}
\usetikzlibrary{arrows}
\usetikzlibrary{automata, positioning}
\usepackage{url}
\usepackage[colorinlistoftodos, color=blue!20]{todonotes}

\usepackage[colorlinks=true,
		linkcolor=MidnightBlue,
		citecolor=Mahogany,
		urlcolor=ForestGreen,
		pdfauthor={Priyanka Dey},
		pdftitle={},
		pdfkeywords={},
		pdfsubject={Technical Report},
		plainpages=false]{hyperref}

\allowdisplaybreaks
\newtheorem{theorem}{Theorem}[section]
\newtheorem{proposition}[theorem]{Proposition}
\newtheorem{lemma}[theorem]{Lemma}

\theoremstyle{definition}

\theoremstyle{remark}
\newtheorem{remark}[theorem]{Remark}

\theoremstyle{assumption}
\newtheorem{assumption}[theorem]{Assumption}

\theoremstyle{observation}

\theoremstyle{problem}

\theoremstyle{fact}

\numberwithin{equation}{section}

\DeclareMathOperator*{\minimize}{minimize}

\DeclareMathOperator{\sbjto}{subject\;to}

\newcommand{\norm}[1]{\left\lVert{#1}\right\rVert} 
\newcommand{\abs}[1]{\left\lvert{#1}\right\rvert} 

\newcommand{\R}{\mathbb{R}} 
\newcommand{\Nz}{\mathbb{N}} 
\newcommand{\N}{\Nz^\star} 
\newcommand{\Let}{\coloneqq}

\newcommand{\st}{x}

\newcommand{\sys}{A}

\newcommand{\match}{M}
\newcommand{\GAB}{G(A, B)}
\newcommand{\bil}{V_{A}^{1}}
\newcommand{\bir}{V_{A}^{2}}
\DeclareMathOperator{\fral}{for\;all\;}
\newcommand{\BAB}{\Gamma(A,B)}
\newcommand{\BA}{\Gamma(A)}

\title[Minimum cost sparsest input-connectivity ]{On minimum cost sparsest input-connectivity for controllability of linear systems}

\author[P. Dey]{Priyanka Dey}
\address{Systems \& Control Engineering, Indian Institute of Technology Bombay, Powai, Mumbai~400076, India.}
\email{dey\_priyanka@sc.iitb.ac.in}

\author[N. Balachandran]{Niranjan Balachandran}
\address{Department of Mathematics, Indian Institute of Technology Bombay, Powai, Mumbai 400076, India.}
\email{niranj@math.iitb.ac.in}

\author[D. Chatterjee]{Debasish Chatterjee}
\address{Systems \& Control Engineering, Indian Institute of Technology Bombay, Powai, Mumbai 400076, India.}
\email{dchatter@iitb.ac.in}
\begin{document}
\keywords{structural controllability, graph theory, matching, flow networks.}

\begin{abstract}
		We deal with algorithmic techniques for minimal cost input-connectivity while maintaining controllability of linear systems. The input matrix is assumed to be constrained in the sense that the set of states that each input (if present) can influence is known a priori, and that each interconnection between an input and a state is associated with a certain cost. In this setting we determine a set of input-connections that lead to the minimum cost and ensures that the resulting system is structurally controllable. We also identify a sparsest set of input-connections with minimum cost while maintaining structural controllability of the system. A large class of systems are identified for which these problems are solvable in polynomial time using efficient algorithms. A 2-approximation solution is presented for the general case. Graph-theoretic tools are employed to tackle the above class of constrained design problems. Illustrative examples are included to demonstrate the efficacy of the techniques developed here. 
	\end{abstract}
	
	\maketitle
	

\section{Introduction}
\label{s:intro}

Dynamical networks arise in a wide-range of application scenarios involving systems such as biological, social, economical, industrial, and transportation systems \cite{ref:OraMoeMur-10, ref:MarKai-12, ref:CurJacPin-09, ref:WatStr-98}, where the states of the systems are updated overtime via its own dynamics. For example, a traffic network \cite{ref:Liu-05} can be modelled as a dynamical system where the load on each road get influenced by other nearby roads and therefore needs to be updated frequently, a social network \cite{ref:OpaPan-09} where the state of each person (describing, e.g., his/her opinion on a particular topic,) gets affected by his/her friends and is updated from time to time, a gene regulatory network \cite{ref:ISte-04} where expressions of proteins are affected by certain specific genes, etc. Of late, the subject of control of such large-scale systems has been attracting considerable attention due to its important academic and practical significance. 

On the one hand, the gigantic sizes of the large-scale systems available today have made the problem of identifying a smallest subset of the inputs to control the systems a very relevant problem. This particular problem is difficult: in fact, it was proved in \cite{ref:AOls-14} that the problem of finding  a smallest set of actuators to ensure a linear system controllable is NP-hard. \footnote{Researchers have also looked at a quantative approach to analyze the ``extent'' of the controllability of a linear system via some energy related metrics such as worst-case control energy, etc. Various tools are employed to find a set of at most \(k\) states to be actuated to optimize these metrics \cite{ref:SumCorLyg-16, ref:PasZamBul-14}.} On the other hand, in many practical situations, it is often necessary to consider the cost of the interconnection connecting an input to a system state; this cost may depend on various factors including the specific functionality of each control, reliability, installation and maintenance, and even environmental conditions such as humidity and temperature. In such cases it is necessary to minimize the number of connections between the inputs and the states as well as the cost of using those connections. We consider several examples: The first is a system of robots modelled as a dynamical system, where each robot communicates with some others over wireless channels to perform a pre-defined task. Assume that a fixed set of external inputs are directly connected to a pre-specified set of states of this system. Each connection has a cost that may depend on various possible factors. In such a setting it is desirable to have relatively a minimal cost set of connections directly connected to a collection of robots rather than squander resources by distributing individual controllers to each robot and indirectly controlling the rest of the system. Such an architecture is especially common in decentralized control and situations where a central controller may be incapable of simultaneously controlling all agents/component subsystems due to, e.g., bandwidth restrictions.

The second example is a multi-agent system, and we consider the situation where the signals are provided to different agents via multiple communication/control towers, and each tower provides signals to multiple agents with cost function that depends on the agent's location and/or characterization. For instance, the cost may be related to the distance between the tower and the agent, or preferential assignment in a heterogeneous multi-agent network where some tower can provide signals to agents of a specific type. Under this setting it is important to optimally select the interconnections between the towers and the agents incurring the least cost while maintaining controllability of the system.

Our third example is that of opinion control in social networks. This topic has attracted considerable attention in recent years due to its widespread applications in psychology, economics, political science, etc., \cite{ref:DasGolMun-14, ref:EteBa-15}. An opinion on social network evolves over time, and the opinion of a particular person is influenced by his/her neighbours (represented in the form of connecting edges). Moreover, such opinions are typically driven by professional journalists, propaganda in the media, campaigns, market agencies, etc., and they act as an external set of inputs to drive the opinions of the people in the social network \cite{ref:DeValGan-16, ref:ArcGri-17}. Steering the opinions of people in social networks is, clearly, of immense importance from a political standpoint, and controllability of such networks is a crucial property. An important particular problem in this setting is that of identifying a minimal set of connections that must be present between the control inputs and the people constituting the social network to drive the opinion in a favourable manner to a desired state.

In the article at hand we focus on \emph{identifying a minimal cost set of connections between the inputs and the system states} to make the system \emph{controllable with a pre-specified input structure}. By input structure we mean that the set of states to which an input is directly connected is known a priori. To this end, we employ structural systems theory to address this problem, where a class of system theoretic problems may be treated by employing only the connections between the system states and inputs. Several interesting and perhaps non-intuitive assertions can be derived via this theory and it is useful especially for systems whose parameters are not exactly known due to various reasons including ageing of system components, structural alterations, etc. Structural analysis of control systems via \emph{structural controllability} was introduced in \cite{ref:Lin-74} for single input linear systems, the ideas were extended to multi-input linear systems in \cite{ref:ShiPea-76}.
Over the past several years a considerable amount of research has been done in this area, see e.g., \cite{ref:May-81, ref:DioComWou-03, ref:LiuSloBar-11, ref:Ols-15, ref:KhaJad-11, ref:DoosKha-13}. 

\emph{Throughout this article we assume that the input matrix and the set of states that each input can influence are known a priori.} We study the techniques of structural system theory to deal with these two different but related problems corresponding to the controllability of a linear system:
\begin{itemize}[label=\(\circ\), leftmargin=*]
\item We assume that each connection between an input and a state has a non-negative cost associated with it. The focus of the first problem is to identify a set of connections with minimum cost while maintaining structural controllability of the system. 
\item The objective of the second problem is to determine a minimal set of connections between the inputs and the states that incurs minimum cost while ensuring the structural controllability of the system.
\end{itemize}
The precise statements of the above problems are given in \S\ref{ss:problems}. We identify mild conditions under which these problems are solvable in polynomial time (in dimension of states and inputs) using efficient algorithmic techniques. A part of this work where the bipartite graph associated with system matrix \(A\) has a perfect matching (see \S\ref{s:Background} for a formal definition) is investigated by us in \cite{ref:DeyBalCha-18a}.

This article unfolds as follows: \S\ref{s:problemformulations} gives the precise statement of the problems under consideration, and related work in this area. \S\ref{s:Background} reviews certain concepts from graph theory that will be needed in this sequel. \S\ref{s:mainresults} provides efficient polynomial time algorithms to obtain solution for the problems stated in \S\ref{ss:problems} under mild assumptions on the system matrix. Moreover, we provide an approximate solution for the general set up. In \S\ref{s:simulations} we demonstrate the effectiveness of our algorithms by providing some illustrative examples. An illustrative benchmark example of the IEEE 118-bus power network is also presented in \S\ref{s:simulations}. 

\section{Problem Formulation}
\label{s:problemformulations}
 The notations employed here are standard: We denote the set of real numbers by \(\R\), the set of non-negative real numbers by \(\R^{+}\), the positive integers by \(\N\), and we let \([n]\Let\lbrace 1,2,\ldots, n\rbrace\) for \(n\in \N\). We denote by \(\abs{X}\) the cardinality of a finite set \(X\). We denote by \(I_n\) the identity matrix of dimension \(n\). If \(A \in \R^{n \times n}\), then \(A_{ij}\) represents the entry located at \(i\)th row and \(j\)th column.
We define \(\textbf{1}\) associated with the \(A_{ij}\) entry as follows:
\begin{equation}
\label{indicator}
\begin{aligned}
&\textbf{1}_{\lbrace A_{ij}\neq 0\rbrace}\Let\begin{cases} 1\;\; \mbox{if}\; A_{ij}\neq 0,\\
 0\;\; \mbox{otherwise}.
\end{cases}
\end{aligned}
\end{equation}
\subsection{Model of the Linear system}\label{ss:Model}
Consider a linear time-invariant system
\begin{equation}
\begin{aligned}
	\label{e:linsys}
	\dot{\st}(t)&= \bar{\sys} \st(t) + \bar{B}u(t),  
\end{aligned}
\end{equation}
where \(\st(t) \in \R^{d}\) are the states and \(u(t)\in \R^{m}\) are the inputs at time \(t\), and \(\bar{\sys} \in \R^{d \times d}\) and \(\bar{B}\in \R^{d\times m}\) are the given state and input matrices respectively. The system \eqref{e:linsys} is completely described by the pair \((\bar{A},\bar{B})\), and we shall interchangeably refer to \eqref{e:linsys} and \((\bar{A},\bar{B})\) in this sequel.  

In our analysis the precise numerical values of the entries of \(\bar{A}\) and \(\bar{B}\) will not matter, but the information about the locations of fixed zeros in \(\bar{A}\) and \(\bar{B}\) will be essential. For any matrix \(R\), its \emph{sparsity matrix} is defined to be matrix of same dimension as \(R\) with each entry either a zero or an independent free parameter, denoted by \(\star\). A \emph{numerical realisation} of \(R\) is obtained by assigning numerical values to the star entries of the sparsity matrix of \(R\).
Let \(A\in \lbrace 0,\star \rbrace^{d \times d}\) and \(B\in \lbrace 0,\star \rbrace^{d \times m}\) represent the sparsity matrices of the system matrix \(\bar{\sys}\) and the input matrix \(\bar{B}\). We say that a pair \((A, B)\) is \emph{structurally controllable} if there exists at least one numerical realization \((A^{\prime}, B^{\prime})\) of \((A, B)\) such that \((A^{\prime}, B^{\prime})\) is controllable.\footnote{It is well-known that if a pair \((A, B)\) is structurally controllable, then \emph{almost all} numerical realizations of \((A, B)\) are controllable \cite{ref:Rei-88}.}

Given a linear time-invariant system \eqref{e:linsys}, a digraph \(G(A,B)\) is associated with it in a natural way: Let \(\mathcal{A}=\lbrace v_1, v_2, \ldots,  v_d\rbrace\) and \(\mathcal{U}=\lbrace u_1, u_2, \ldots, u_m\rbrace\) be the state and the input vertices corresponding to the states \(\st(t)\in \R^d\) and the inputs \(u(t) \in \R^m\), respectively, of the system \eqref{e:linsys}. Let \( E_{A}=\lbrace (v_j, v_i)\,|\, {A}_{ij}\neq 0\rbrace\) and \(E_{B}=\lbrace (u_j, v_i)\,|\, B_{ij}\neq 0\rbrace\). We define a digraph \(G(A,B)=(V, E)\), where \(V=\mathcal{A}\, \sqcup \,\mathcal{U}\), \(E= E_{A} \sqcup E_{B}\), and \(\sqcup\) represents the disjoint union. The edges in \(E_A\) are referred as \emph{state-connections}. The edges in \(E_B\) are referred as \emph{input-connections} in this sequel. Sometimes we shall need the digraph \(G(A)= (\mathcal{A}, E_{A})\) with vertex set \(\mathcal{A}\) and edge set \(E_{A}\) considering the edges between only the state vertices.

\subsection{Problem Statement}
\label{ss:problems}
Before formally stating the two optimisation problems, we introduce the various norms needed in this sequel. 

For a matrix \(N \in \lbrace 0, \star \rbrace^{n \times k}\) (where \(n, k \in \N\))
\begin{itemize}[label=\(\circ\), leftmargin=*]
\item \(\norm{N}_0\) denote the number of non-zero entries in the matrix \(N\).
\item Let each non-zero entry in \(N\) is associated with a non-negative cost. For instance, if \(N_{ij}\neq 0\) assume that \(w_{ij}\geq 0\) be the cost corresponding to it. We define
\[
\norm{N}_{w}\Let \sum_{i=1}^{n}\sum_{j=1}^{k} w_{ij} \textbf{1}_{\lbrace N_{ij}\neq 0\rbrace}, 
\]
where \(\textbf{1}\) is defined in \eqref{indicator}. 
\end{itemize}

\emph{Given \(A \in \lbrace 0, \star \rbrace^{d \times d}\) and \(B\in \lbrace 0, \star \rbrace^{d \times m}\), throughout we assume that the pair \((A,B)\) is structurally controllable.} For a matrix \(B \in \lbrace 0, \star \rbrace^{d \times m}\), let the collection of locations of fixed zeros of \(B\) be \(Z(B)\), i.e., \(Z(B)\Let \big\lbrace (i,j)\;|\; B_{ij}=0\big\rbrace\). Define 
\[\mathcal{K}\Let\Big\lbrace B^{\prime}\,\Big|\, Z(B)\subset Z(B^{\prime}),\, (A, B^{\prime})\,\text{is structurally controllable}\Big\rbrace.
\]
Since \((A,B)\) is structurally controllable by assumption, \(\mathcal{K}\) is always non-empty. 

We deal with the following two optimization problems:\\
\noindent \textbf{Problems:} Let \(A \in \lbrace 0, \star \rbrace^{d \times d}\) and \(B\in \lbrace 0, \star \rbrace^{d \times m}\) be given such that \((A,B)\) is structurally controllable:
\begin{itemize}[label=\(\circ\), leftmargin=*]
\item Each non-zero entry in \(B\) is associated with a non-negative cost. Let \(w_{ij}\geq 0\) denote the cost of using the input-connection connecting the input vertex \(u_j\) to the state vertex \(v_i\) in \(G(A,B)\). Determine an input matrix \(B^{*}\) that solves the problem:
 \begin{equation*}
\begin{aligned}
\minimize_{B^{\prime}\in \mathcal{K}}\;\; \norm{B^{\prime}}_w,
\end{aligned}
\tag{\(\mathcal P_1\)}
\label{P3}
\end{equation*}
\begin{itemize}[label=\(\rhd\), leftmargin=*]
\item The following is a special case of \eqref{P3}: Assume that a fixed and non-zero cost is assigned to all the input-connections in \(G(A, B)\). Determine an input matrix \(B^{*}\) that solves
\begin{equation*}
\begin{aligned}
\minimize_{B^{\prime}\in \mathcal{K}} && \norm{B^{\prime}}_0.
 \end{aligned}
 \tag{\(\mathcal P_1^{\prime}\)}
 \label{P1}
\end{equation*}
\end{itemize}
\item Let \(w_{ij}\geq 0\) denote the cost of using the input-connection connecting the input vertex \(u_j\) to the state vertex \(v_i\) in \(G(A,B)\). Determine an input matrix \(B^{*}\) that solves the following optimisation problem:
\begin{equation*}
\begin{aligned}
&\minimize_{B^{\prime}\in \mathcal{K}}\;\; \norm{B^{\prime}}_w\\
& \sbjto \; \norm{B^{\prime}}_0 \leq \norm{B^{\prime\prime}}_0 \;\;\fral B^{\prime\prime} \in \mathcal{K}. \\ 
  \end{aligned}
\tag{\(\mathcal P_2\)}
\label{P2}
\end{equation*} 
\end{itemize}
\noindent where \(\norm{B^{\prime}}_w= \sum_{i=1}^{d}\sum_{j=1}^{m}w_{ij} \textbf{1}_{\lbrace B_{ij}^{\prime}\neq 0\rbrace}\).

\begin{remark}
\label{brute}
As part of our premise, we assume throughout that the given pair \((A,B)\) is structurally controllable. Clearly, finding brute-force solutions to Problems \eqref{P3}-\eqref{P2} requires checking all possible input matrices \(B^{\prime}\in \mathcal{K}\), which is quite an impossible combinatorial problem for even moderately sized pairs \((A, B)\). For example, consider  an input matrix \(B \in \lbrace 0, \star \rbrace^{d \times m}\) containing \(n\) non-zero entries. To identify a solution to \eqref{P3}-\eqref{P2} requires testing all possible combinations of the subsets of the  \(n\) entries. Therefore, the number of computations needed is exponential. However, the algorithms discussed here identify solutions of these problems efficiently in polynomial time complexity without using brute-force techniques.
\end{remark} 
\vspace*{-1cm}
\noindent \textbf{Related work:} We restrict our attention to only those articles that are closely related to our work; structural analysis for various other problems can be found in \cite{ref:LiuBarAlb-16} and the references therein. Let \(d\) be the number of states in the system. The problem of finding the \emph{minimum number of inputs} required to guarantee structural controllability was considered in \cite{ref:ComDioWou-02}, \cite{ref:LiuSloBar-11}, but these results do \emph{not} solve \eqref{P3} and \eqref{P2} since the former problem gives only a lower bound on the minimum number of inputs present in our solutions and gives no information about the set of \emph{input-connections} (of minimal cost) required when we have a pre-specified input structure. The problem of \emph{selecting the fewest states to be influenced} to achieve structural controllability, when the given input matrix is unconstrained, is considered in the single-input case in \cite{ref:ComDio-15} and for the multi-input case in \cite{ref:PeqKarAgui-16}, and the same problem is dealt in \cite{ref:Ols-15} to give a solution in lower computational complexity. This particular problem is a special case of \eqref{P3}, namely, when all the input-connections are permissible between the given inputs and states with a uniform non-zero cost assigned to each input-connection. Since there is no constraint on the available set of input-connections in \cite{ref:PeqKarAgui-16}, their proposed solutions (in general) are not useful to solve \eqref{P3} even subproblem \eqref{P1}. The article\cite{ref:PequKarAgu-16} considered the \emph{minimum cost input design problem} where the objective is to find an input matrix which makes the system structurally controllable and incurs minimum cost when each \emph{state} is associated with a certain cost that is independent of the input performing the task. In constrast, we aim to minimize the cost of actuating input-connections, where the cost depends on the input-to-state pair; our problem is more general while satisfying the distributed controllability constraints. In fact, the setting of \cite{ref:PequKarAgu-16} can be reduced to the setting of \eqref{P3} by assuming that the given input matrix is an unconstrained \(d\times d\) matrix with identical cost given to all the \(d\) input-connections associated to a state.

Apart from these prior investigations, when the given input matrix is \emph{constrained}, then the problem of selecting an input set of minimum cardinality to guarantee that the resulting system is structurally controllable is considered in \cite{ref:ComDio13a} and shown to be NP-hard in \cite{ref:PeqKarAgu-15}. This problem is extended to the case where a non-negative cost is assigned to each given input in  \cite{ref:PeqKarPap-15a}, \cite{ref:MooChaBel-18}, known as \emph{minimum cost constrained input selection problem} (minCCIS). This problem is reduced to the minimum cost fixed flow problem and a \emph{polynomial time \(\Delta\)-approximation solution}, where \(1\leq \Delta\leq d+1\) is provided in \cite{ref:MooChaBel-18}. The problems tackled in this article are quite distinct from the minCCIS: The latter aims at identifying a minimal cost input set whereas our problems concern is \emph{identifying a minimal cost set of the input-connections from the already deployed set of input-connections in the system.} Remarks  \ref{r:minCISP1}, \ref{differencestrong}, and  \ref{r:differenceperfect} contain further technical discussions.

\section{Background}
\label{s:Background}

Given \(G(A,B)\), a state vertex \(v_i \in \mathcal{A}\) is said to be \emph{accessible} if there exists a directed path\footnote{A sequence of edges \(\lbrace(v_1, v_2), (v_2, v_3), \ldots, (v_{k-1}, v_k)\rbrace\) with all vertices distinct is called a \emph{directed path} from \(v_1\) to \(v_k\).} from some input \(u_j\) to \(v_i\); otherwise, it is \emph{inaccessible}. A \emph{cycle} is a directed path where the initial vertex \(v_1\) coincides with the end vertex \(v_k\). The digraph \(G(A)\) is \emph{strongly connected} if for each ordered pair of vertices \((v_i, v_j)\), there exists a directed path from \(v_i\) to \(v_j\). A \emph{strongly connected component} (SCC) of \(G(A)\), usually denoted by \(\mathcal{S}\), is a maximal strongly connected subgraph of \(G(A)\). An SCC \(\mathcal{S}\) in the digraph \(G(A)\) is said to be \emph{source strongly connected component} (SSCC) if there is no directed edge from the vertices of other SCCs into any vertex of \(\mathcal{S}\). As a consequence, all the states of \(G(A)\) are accessible if and only if all the SSCCs are accessible. In addition to the accessibility condition the digraph \(G(A,B)\) should also satisfy a \emph{no-dilation} condition. \footnote{The digraph \(G(A,B)\) said to have dilation if there exists a set of state vertices \(T\) whose neighbourhood set \(N^{-}(T)\)(where a vertex \(v\in N^{-}(T)\), if there exists a directed edge from \(v\) to a vertex in \(T\)) has fewer vertices than \(T\). Please refer to \cite{ref:Lin-74, ref:LiuSloBar-11} for more details.} A fundamental connection between system theoretic property of structural controllability and certain structural properties of \(G(A,B)\) is given by:
\begin{theorem}\cite[Theorem 1, p.~207]{ref:Lin-74} 
\label{t:lin} The pair \((A,B)\) is structurally controllable if and only if the associated digraph \(G(A,B)\) derived from \eqref{e:linsys} has every state vertex \(v_i\in \mathcal{A}\) accessible and  is free of dilations.
\end{theorem}

The digraph \(G(A, B)\) derived above from \eqref{e:linsys} can also be represented by an undirected bipartite graph in the following standard fashion:
\(\BAB \Let((\bil \sqcup V_{B}), \bir, \mathcal{E}_A \sqcup \mathcal{E}_B)\), where \(\bil\Let{\lbrace v_1^1, v_2^1, \ldots, v_d^1 \rbrace} \), \(\bir\Let{\lbrace v_1^2, v_2^2, \ldots, v_d^2 \rbrace} \), and \(V_B\Let \lbrace u_1, u_2, \ldots, u_m\rbrace\). \( \mathcal{E}_A = \lbrace (v_j^1, v_i^2)\,|\, A_{ij}\neq 0\rbrace \) and \( \mathcal{E}_B = \lbrace (u_j, v_i^2)\,|\, B_{ij}\neq 0 \rbrace\). In the similar manner, we can define the undirected bipartite graph \(\BA\Let((\bil , \bir), \mathcal{E}_A)\) which is the induced subgraph\footnote{An induced subgraph of a graph is another graph, formed from a subset of the vertices of the graph and all of the edges connecting pairs of vertices in that subset.} on \(\bil \sqcup \bir\) of \(\BAB\). We shall need a few more definitions in the context of graph \(\BA\). A \emph{matching} \(\match\) in \(\BA\) is a subset of edges that do not share vertices. A \emph{maximum matching} \(\match\) in \(\BA\) is a matching with maximum number of edges. An edge \(e\) is said to be \emph{matched} if \(e\in M\). A vertex is said to be \emph{matched/saturated} if it belongs to an edge in the matching \(\match\); otherwise, it is \emph{unmatched.} A matching \(\match\) in \(\BA\) is said to be \emph{perfect} if all the vertices in \(\BA\) are matched. We say that \(\BAB\) has a \emph{system of distinct representatives (SDR)} if there exists a matching \(M\) that saturates all the vertices of \(\bir\) in \(\BAB\). The presence of dilation in \(G(A,B)\) can be easily checked by using a matching condition that relates \(\BAB\) and the no-dilation condition.
\begin{proposition}\cite[Theorem 2]{ref:Ols-15}
\label{p:Aols}
A digraph \(G(A,B)\) has no dilation if and only if there exists an SDR in the bipartite graph \(\BAB\).
\end{proposition}  
\begin{remark}
\label{r:structure}
Given \(G(A)=(\mathcal{A}, E_A)\), the SSCCs can be determined in \(O(\abs{\mathcal{A}}+\abs{E_A})\) computations  \cite{ref:CorLeiRivSte-09}. We know that \(\abs{\mathcal{A}}=d\) and \(\abs{E_A}=O(d^2)\). The procedure for finding the SSCCs involves \(O(d^2)\) computations and checking for existence of an SDR in \(\BAB\) involves \(O(d^{2.5})\) computations. Thus, structural controllability of a pair \((A,B)\) can be accurately checked in  \(O(d^{2.5})\) computations \cite{ref:Die-20}.
\end{remark}

Given \(\BAB\), a cost function \(c: \mathcal{E}_A \sqcup \mathcal{E}_B \to \R^{+}\) assigns non-negative costs to the edges of the bipartite graph \(\BAB\), represented by \((\BAB ;c)\). A minimum cost maximum matching(MCMM) finds a maximum matching \(M^{*}\) such that \(\sum_{e \in M^{*}} c(e) \leq \sum_{e \in \bar{M}} c(e)\), where \(\bar{M}\) is any maximum matching in \(\BAB\) \cite{ref:Die-20}.

\section{Main Results}
\label{s:mainresults}

\subsection{Systems with perfect matching} \label{ss:perfect matching}

\begin{assumption}
\label{a:perfectmatching}
We stipulate that the system matrix \(A\) is such that the bipartite graph \(\BA\) has a perfect matching.
\end{assumption}
This is indeed a reasonable assumption since a large class of systems, for instance systems including epidemic dynamics, power grids, multi-agent systems,  etc., \cite{ref:IliXieKhaMou-10, ref:JafAjoAgh-11, ref:NewBarDun-06} exhibit this feature. Assumption \ref{a:perfectmatching} ensures that the bipartite graph \(\BA\) has a perfect matching \(M\). Since \(\BA\) is an induced subgraph on \(\bil \sqcup \bir\) in \(\BAB\),  matching \(M\) saturates all the vertices of \(\bir\) in \(\BAB\).  Thus, the bipartite graph \(\BAB\) has an SDR. Proposition \ref{p:Aols} implies that \(G(A,B)\) has no-dilation even in the absence of any input vertex. By Theorem \ref{t:lin}, it suffices to ensure the accessibility criterion for structural controllability of the pair \((A, B)\). This implies that the pair \((A,B)\) is structurally controllable if and only if all the SSCCs of \(G(A)\) are accessible from some input vertex in \(G(A,B)\). The following Lemma is central to the development of our results:
\begin{lemma}
\label{l:sparseBperfectmatching}
Suppose that Assumption \ref{a:perfectmatching} holds. Consider a structurally controllable pair \((A,B)\), and let \(q\) denote the number of SSCCs in \(G(A)\). If \(B^{*}\) solves \eqref{P1}, then \(\norm{B^{*}}_0= q\).
\end{lemma}

We move on to \eqref{P2}. Recall from \S\ref{ss:problems} that \(w_{ij}\) is the cost associated with the input-connection from input vertex \(u_j\) to state vertex \(v_i\). Let \(w_{\max}\) denote the maximum cost assigned to an input-connection (corresponding to a non-zero entry in \(B\),) among all the input-connections present in \(G(A, B)\). In our setting for solving \eqref{P2}, we impose the condition that if an input vertex \(u_k\) does not have an input-connection to a state vertex \(v_{\ell}\) (determined by the given input matrix \(B\)), then \(w_{\ell k}= \infty\) (for practical purposes \(w_{\ell k}\) is taken to be \(w_{\max}+1\)). We provide the following Algorithm \ref{AlgorithmProblem2} to obtain a solution of \eqref{P2}.

\begin{algorithm}[h]
\KwIn{\(A \in \lbrace 0, \star \rbrace^{d \times d}\), \(B \in \lbrace 0, \star \rbrace^{d \times m}\)}
	\KwOut{The input matrix \(B^{*}\)}
	
	\nl Determine the SSCCs \(\lbrace\mathcal{S}_j\rbrace_{j=1}^{q}\) 
	
	\nl \(L \leftarrow \emptyset\) 
	
	\nl for each \(\mathcal{S}_j\) \textbf{do}
	
	\begin{enumerate}[leftmargin=*]
	\item for each state vertex \(v_i \in \mathcal{S}_j\), choose the smallest cost \(w_{ik}\) among \(\lbrace w_{i1}, w_{i2},\ldots, w_{im}\rbrace\).
	\item choose a state vertex of least cost, say \(v_\ell\), with cost \(w_{\ell k}\) among all the state vertices in \(\mathcal{S}_j\).
	\end{enumerate}
	
	\nl \(L \leftarrow L \cup (u_k, v_{\ell})\)
	
	\nl end \textbf{for}
	
	\nl Define:
	 \begin{equation*}
\label{e:SCSProblem2}
B^{*}_{\ell k}\leftarrow\begin{cases} \star \, &\mbox{if}\,\, e=(u_k, v_{\ell}) \in L,\\
0\, &\mbox{otherwise}.
\end{cases}   
\end{equation*}
	\caption{Algorithm to solve \eqref{P2}}
\label{AlgorithmProblem2}
\end{algorithm}

The structural controllability of the given pair \((A, B)\) ensures that no input-connection corresponding to cost \(\infty\) (\(w_{\max}+1\)) is selected by the algorithm.
\begin{theorem}
\label{t:problemP2}
Let \((A, B)\) be a linear system and suppose that Assumption \ref{a:perfectmatching} holds. The procedure outlined in Algorithm \ref{AlgorithmProblem2} yields a matrix \(B^{*}\) such that \((A,B^{*})\) is structurally controllable and solves \eqref{P2}. 
\end{theorem}
Finding the SSCCs involves \(O(d^2)\) computations \cite{ref:CorLeiRivSte-09} (Step 1). The overall complexity of Step 3-5 is linear in number of edges in \(E_B\). Since \(\abs{E_B}\leq dm\) the complexity of Algorithm \ref{AlgorithmProblem2} to identify a solution of \eqref{P2} is \(O(dm)\).

By setting unit cost to all the input-connection in \(G(A, B)\), we can obtain a solution of \eqref{P1} by using Algorithm \ref{AlgorithmProblem2}. 
The strategy designed for identifying a solution to \eqref{P2} namely, Algorithm \ref{AlgorithmProblem2}, also provides a solution to \eqref{P3}. 
\begin{proposition}
\label{P:problem3assum1}
Let \((A, B)\) be a linear system and suppose that Assumption \ref{a:perfectmatching} holds. Then \(B^{*}\) obtained in Algorithm \ref{AlgorithmProblem2} also solves \eqref{P3}.
\end{proposition}

\begin{remark}
\label{r:minCISP1}
Generally, the vertex-variant of a problem is difficult to solve as compared to the edge-variant associated with it. For example, finding an independent set (a set of non-adjacent vertices) of maximum cardinality in an undirected graph is an NP-hard problem. However, the problem finding a maximum matching (a set of non-adjacent edges) in an undirected graph admits many polynomial time algorithms to compute it optimally. In a similar manner, a different but related problem is the minCCIS discussed in \S\ref{s:problemformulations}. The authors of \cite{ref:PeqKarAgu-15} observed that the minCCIS is NP-hard when the bipartite graph \(\BA\) associated with the system matrix \(A\) has a perfect matching. However, \emph{we demonstrate in this sequel that it is possible to select a minimal cost input matrix \(B^{*}\) in polynomial time (when the input matrix \(B\) is known a priori.}
\end{remark}
\subsection{Strongly connected systems}
\label{ss:stronglyconnected}
\begin{assumption}
\label{a:stronglyconnected}
We stipulate that the system matrix \(A\) is such that the digraph \(G(A)\) is strongly connected.
\end{assumption}
There exists a variety of interconnected dynamical system \cite{ref:Ant-13}, dynamics based on consensus-like protocol \cite{ref:DimKArMouRabSca-10}, \cite{ref:JadLinMor-03} where strong connectivity is essentially ensured by network connectivity. The preceding condition assures that all the state vertices in \(G(A)\) are accessible by using only one input-connection connecting an input vertex to a state vertex in \(G(A,B)\). If the system matrix \(A\) is such that both Assumption \ref{a:perfectmatching} and Assumption \ref{a:stronglyconnected} are satisfied, then Problems \eqref{P3}-\eqref{P2} can be solved in a straightforward manner: Only one input-connection is enough to satisfy the conditions given in Theorem \ref{t:lin}, i.e., to ensure structural controllability of pair \((A,B)\). Thus, \(\norm{B^{*}}_0=1\). Also, we may pick an input-connection with the least cost among the costs assigned to all the input-connections. Thus, \emph{we further assume that the system matrix \(A\) is such that the bipartite graph \(\BA\) does not have perfect matching.} As a consequence, by Theorem \ref{t:lin}, only the no-dilation criterion has to be satisfied to ensure structural controllability of the pair \((A,B)\). 

We start with an algorithm to solve \eqref{P3}:
\begin{algorithm}
\KwIn{\(A \in \lbrace 0, \star \rbrace^{d \times d}\) and \(B \in \lbrace 0, \star \rbrace^{d \times m}\)}
	\KwOut{The input matrix \(B^{*}\)}
	
	\nl construct \(\BAB=(\bil \sqcup V_{B}, \bir, \mathcal{E}_A \sqcup \mathcal{E}_B)\).
	
	\nl for each edge \(e \in \BAB\) define cost \(c\)
	
	\nl define:
	\begin{equation}
\label{e:weightProblem3assum4.2}
	c(e) \leftarrow \begin{cases} 0 \, & \mbox{for}\,\, e=(v_r^1,v_k^2) \in \mathcal{E}_A,\\
	w_{kj}\, & \mbox{for}\,\, e=(u_j,v_k^2) \in \mathcal{E}_B.
	\end{cases}
	\end{equation}
	
	\nl find a minimum cost maximum matching(MCMM)\cite{ref:JMun-57} in \((\BAB ;c)\), say \(M^{*}\)
	
	\nl define:
	 \begin{equation*}
\label{e:solutionProb2}
B^{*}_{\ell k}\leftarrow\begin{cases} \star \, &\mbox{if}\,\, e=(u_k, v_{\ell}^2) \in M^{*} \cap \mathcal{E}_B,\\
0\, &\mbox{otherwise}.
\end{cases}   
\end{equation*}
\caption{Algorithm to solve \eqref{P3}}
\label{algorithmstronglyprob2}
\end{algorithm}
\begin{theorem}
\label{t:assumption4.2Prob2}
Consider a linear system \((A,B)\) and Suppose Assumption \ref{a:stronglyconnected} holds. The procedure outlined in Algorithm \ref{algorithmstronglyprob2} ensures that the input matrix \(B^{*}\) obtained from it is such that \((A,B^{*})\) is structurally controllable and solves \eqref{P3}. In addition, the complexity of the algorithm is \(O((d+m)^3)\).
\end{theorem}

Algorithm \ref{algorithmstronglyprob2} is also utilized to provide a solution \(B^{*}\) for \eqref{P1} when Assumption \ref{a:stronglyconnected} holds. 
 If all the input-connections have uniform non-zero cost, then \eqref{P1} is a special case of \eqref{P3}. We define the following cost function \(c\) from \(e \in \BAB\) for Algorithm \ref{algorithmstronglyprob2} to solve \eqref{P1}.  
\begin{equation}
\label{e:weightProblem1assum4.2}
	c(e) \leftarrow \begin{cases} 0 \, & \mbox{for}\,\, e \in \mathcal{E}_A,\\
	1\, & \mbox{for}\,\, e \in \mathcal{E}_B.
	\end{cases}
	\end{equation}
 
Now we move to \eqref{P2}. To address this problem, under Assumption \ref{a:stronglyconnected}, we use the following lemma.
\begin{lemma}
\label{l:solveP2P3}
Consider the linear system \((A, B)\), and suppose Assumption \ref{a:stronglyconnected} holds. If all the input-connections corresponding to the non-zero entries of \(B\) have positive cost then \(B^*\) obtained by Algorithm \ref{algorithmstronglyprob2} solves \eqref{P2}.
\end{lemma}

Lemma \ref{l:solveP2P3} implies that Algorithm \ref{algorithmstronglyprob2} provides a solution to both  \eqref{P2} as well as \eqref{P3} whenever the input-connections in \(\GAB\) have positive cost. However, if any input-connection in \(\GAB\) has zero cost then a solution obtained by Algorithm \ref{algorithmstronglyprob2} may not have least number of non-zero entries. In this situation, we modify the cost function \(c\) given by \eqref{e:weightProblem3assum4.2} in Algorithm \ref{algorithmstronglyprob2} by 
	\begin{equation}
\label{e:weightProblassum4.2}
	c_1(e) \leftarrow \begin{cases} 0 \, & \mbox{for}\,\, e=(v_r^1,v_k^2) \in \mathcal{E}_A,\\
	w_{kj}+1\, & \mbox{for}\,\, e=(u_j,v_k^2) \in \mathcal{E}_B.
	\end{cases}
	\end{equation}
	
\begin{lemma}
\label{o:c1andc2}
Consider a linear system \((A,B)\) and let \(\BAB=(\bil \sqcup V_{B}, \bir, \mathcal{E}_A \sqcup \mathcal{E}_B)\) be its associated bipartite graph such that it admits an SDR for \(\bir\). Let \(c\) and \(c_1\) defined in \eqref{e:weightProblem3assum4.2} and \eqref{e:weightProblassum4.2} be the costs assigned to the edges of  \(\BAB\). If \(M^{*}\) is an MCMM under cost function \(c_1\) then it is also an MCMM under cost function \(c\).
\end{lemma}	
Let \(B^{*}\) be the input matrix obtained by Algorithm \ref{algorithmstronglyprob2} under cost function \(c_1\) defined in \eqref{e:weightProblassum4.2}. Clearly, \(B^{*}\in \mathcal{K}\). Let \(M^{*}\) be an MCMM associated with \(B^{*}\). The cost of \(B^{*}\), i.e., \(\norm{B^{*}}_w\) is equal to the sum of the cost of the edges of \(M^{*}\cap \mathcal{E}_B\) in \(\BAB\). By Lemma \ref{l:solveP2P3} and Lemma \ref{o:c1andc2}, it follows that \(B^{*}\) solves \eqref{P2}.
\begin{remark}
\label{differencestrong}
Suppose Assumption \ref{a:stronglyconnected} holds. The minCCIS is no longer NP-hard for this class as shown in \cite{ref:PeqKarPap-15a}. It suffices to ensure the existence of an SDR in \(\BAB\) when \(G(A)\) is strongly connected and \(\BA\) does not have perfect matching for structural controllability of a pair \((A, B)\). As a consequence, the aim of both the minCCIS and \eqref{P3} is to find a matching \(M\) in \(\BAB\) such that all the vertices of \(\bir\) are saturated by \(M\) while minimizing the given objective function. However, since the objective functions to minimize for  minCCIS and \eqref{P3} are different, an optimal matching \(M\) in \(\BAB\) for the minCCIS may not be optimal for \eqref{P3}, and thus not sufficient to arrive at a solution of \eqref{P3}.
\end{remark}
\begin{remark}
\label{r:differenceperfect}
Suppose Assumption \ref{a:perfectmatching} holds. The article \cite{ref:PeqKarAgui-16} provides insights about the number of input-connections that must be present in the solution of \eqref{P3} and \eqref{P2} obtained by Algorithm \ref{AlgorithmProblem2}, i.e., equal to the number of SSCCs. However, such informations are useful but not sufficient to fully determine the set of input-connections that solves \eqref{P3} or \eqref{P2} since neither it considers the pre-specified input-structure nor the cost associated with each input-connections in our setting. Similar kinds of arguments holds for Assumption \ref{a:stronglyconnected} also.
\end{remark}
\subsubsection{Directed tree}
\label{s:directed tree}
\begin{assumption}
\label{a:directed tree}
We stipulate that the system matrix \(A\) is such that the digraph \(G(A)\) is directed rooted tree at \(v_i\in \mathcal{A}\).\footnote{A digraph is said to be \emph{acyclic} if it contains no cycle. A directed tree rooted at \(x\) is an acyclic graph where \(x\) has a directed path to every vertex in the tree.}
\end{assumption}
Assumption \ref{a:directed tree} ensures that \(G(A)\) has only one SSCC, i.e., \(G(A)\) is rooted at \(v_i\) then \(\mathcal{S}_1=\lbrace v_i \rbrace\). Also, \(v_i\) is unmatched by any matching in \(\BA\) since in-degree of \(v_i\), \(\deg^{-}(v_i)=0\) in \(G(A)\). As a consequence of these two observations, it is easy to see that a pair \((A, B)\) is structurally controllable if and only if there exists an SDR for \(\bir\) in \(\BAB\) when Assumption \ref{a:directed tree} holds. In this case, Algorithm \ref{algorithmstronglyprob2} can be used to solve \eqref{P3} when Assumption \ref{a:directed tree} holds. Note that it is not necessary that \(G(A)\) has to be strongly connected for the application of Algorithm \ref{algorithmstronglyprob2}.

\begin{proposition}
\label{p:directedtreeP3}
Consider a linear system \((A,B)\) and suppose Assumption \ref{a:directed tree} holds. Algorithm \ref{algorithmstronglyprob2} provides an input matrix \(B^{*}\) that solves \eqref{P3}.
\end{proposition}  
To this end, the procedure and results involved to solve \eqref{P2} under Assumption \ref{a:directed tree} is along the same lines as the procedure given to solve \eqref{P2} when Assumption \ref{a:stronglyconnected} holds. Hence, the proofs are omitted.
\begin{remark}
The set of state vertices \(F\subset \mathcal{A}\) is \emph{forbidden} if no input is allowed to be directly connected to any vertex in \(F\).\footnote{The problem of finding an input matrix when a forbidden set is given in dealt in \cite{ref:Ols-15} in the unconstrained setting, i.e., when the given input matrix is \(I_{d}\).} If a non-empty forbidden set is present, then the premise of Problems \eqref{P3}-\eqref{P2} is altered by removing the input-connections that connects an input vertex to a state vertex \(v_i \in F\). We obtain a new input matrix from \(B\), say \(\tilde{B}\). Verification of whether the pair \((A, \tilde{B})\) is structurally controllable or not can be done in polynomial time (see Remark \ref{r:structure}), and the solutions of \eqref{P3}-\eqref{P2} are obtained by using the same techniques as discussed in this article.    
\end{remark}
\subsection{General case}
Recently, the authors of \cite{ref:ZhaZho-17} considered the problem of identifying a sub-collection of state-connections and input-connections, from the available set of state and input-connections, that leads to a minimum cost and guarantees that the resulting system is structurally controllable. It was observed that this problem is NP-hard by showing that the \emph{Hamiltonian path problem}\footnote{A Hamitonian path is a directed path which visits every vertex exactly once. It is well-known that determining whether a  Hamiltonian path exists in a graph is NP-hard \cite{ref:CorLeiRivSte-09}.} is polynomially reducible to an instance of this problem. If zero costs are assigned to the state-connections then this problem reduces to the problem of identifying a sub-collection of input-connections from the available set of input-connections that incurs minimum cost and ensures structural controllability of the resulting system, which is \eqref{P3} precisely. This reduction does not show that \eqref{P3} is also NP-hard. The hardness of \eqref{P3} in the general setup is considered as an open problem and needs further investigation. We do however provide an algorithm which guarantees a \(2\)-approximate solution to \eqref{P3} in the general set-up.  
\begin{algorithm}
\KwIn{\(A \in \lbrace 0, \star \rbrace^{d \times d}\) and \(B \in \lbrace 0, \star \rbrace^{d \times m}\)}
	\KwOut{The input matrix \(B^{*}\)}
	
	\nl Use Algorithm \ref{AlgorithmProblem2} to obtain \(B^{\prime}\).
	
	\nl Use Algorithm \ref{algorithmstronglyprob2} to obtain \(\widetilde{B}\).
	
     \nl Define an input matrix \(B^{\prime\prime}\Let B^{\prime}\)	
	
	\nl Let \(\lbrace\mathcal{S}_j\rbrace_{j=1}^{r}\) be the SSCCs for which there exists a \(v_{\ell} \in \mathcal{S}_j\) s.t. \(\widetilde{B}_{\ell k}=\star\), where \(r \leq q\) and \(k \in [m]\). 
	
	\nl for \(j=1, \ldots, r\)\\
	    update \(B^{\prime\prime}\) by letting \(B^{\prime\prime}_{\ell k}=0\) for all \(v_{\ell} \in \mathcal{S}_j\) and \(k\in[m]\).\\
	    end
	   
	 \nl \(B^{*}=B^{\prime\prime}\sqcup \widetilde{B}\).   
\caption{Finding an approximate solution to \eqref{P3}}
\label{algorithmapproximate}
\end{algorithm}

In Algorithm \ref{algorithmapproximate}, Algorithm \ref{AlgorithmProblem2} is utilised to ensure input-reachability of the state vertices from inputs and Algorithm \ref{algorithmstronglyprob2} is used to obtain an input matrix that removes dilation. Step \(3\) defines a new input matrix \(B^{\prime\prime}\) which is same as \(B^{\prime}\). Step \(4\) collects those SSCCs, say \(\lbrace\mathcal{S}_j\rbrace_{j=1}^{r}\) (where \(r\leq q\)), of \(G(A)\) which have at least one state vertex \(v_{\ell} \in \mathcal{S}_j\) that has an input-connection from some input \(u_k\) in \(\widetilde{B}\), i.e., \(\widetilde{B}_{\ell k}=\star\). This confirms that the SSCCs \(\lbrace\mathcal{S}_j\rbrace_{j=1}^{r}\) are accessible by using the input-connections corresponding to \(\widetilde{B}\). Consequently, Step \(5\) removes those input-connections corresponding to \(B^{\prime\prime}\) needed to make \(\lbrace\mathcal{S}_j\rbrace_{j=1}^{r}\) accessible. The rest of the SSCCs \(\lbrace\mathcal{S}_j\rbrace_{j=r+1}^{q}\) remain accessible by using the input-connections corresponding to the updated \(B^{\prime\prime}\)(obtained after Step \(5\)).

\begin{theorem}
\label{t:approximate}
Consider a linear system \((A, B)\). Algorithm \ref{algorithmapproximate} yields a 2-approximation solution to \eqref{P3} with complexity \(O((d+m)^3)\), where \(d\) and \(m\) are the number of states and inputs in the system.
\end{theorem}

\begin{remark}
By duality \cite[p.~53]{ref:Kai-80} between controllability and observability in linear time-invariant systems, structural controllability of \((A,B)\) is equivalent to structural observability of \((A^{\top},B^{\top})\). Thus, the above results can be directly extended to the corresponding observability problems.
\end{remark}
\section{Illustrative examples}
\label{s:simulations}

\noindent \textit{Example 1}: Let the structures of the system and input matrices be
\[
A=\begin{bmatrix}
0 & \star & 0 &  0 & 0 & 0 & 0 & 0 & 0 & 0\\
0 & \star & \star &  0 & 0 & 0 & 0 & 0 & 0 & 0\\
\star & 0 & 0 &  0 & 0 & 0 & 0 & 0 & 0 & 0\\
0 & 0 & \star &  \star & 0 & \star & 0 & 0 & 0 & 0\\
0 & 0 & 0 &  \star & 0 & 0 & 0 & 0 & \star & 0\\
0 & 0 & 0 &  0 & \star & 0 & 0 & \star & 0 & 0\\
0 & 0 & 0 &  0 & 0 & 0 & 0 & \star & 0 & 0\\
0 & 0 & 0 &  0 & 0 & 0 & \star & 0 & 0 & 0\\
0 & 0 & 0 &  0 & 0 & 0 & 0 & 0 & 0 & \star\\
0 & 0 & 0 &  0 & 0 & 0 & 0 & 0 & \star & 0\\
\end{bmatrix}
,\;
B=\begin{bmatrix}
\star & 0  & 0 \\
0 & 0 & 0 \\
\star & 0  & 0 \\
\star & 0  & 0 \\
0 & 0 & 0 \\
0 & \star  & 0 \\
0 & \star  & 0 \\
0 & 0  & \star \\
0 & 0  & 0 \\
0 & 0  & \star \\
\end{bmatrix}.
\]

\noindent The digraph \(G(A, B)\) associated with the pair \((A, B)\) is shown in Fig.~\ref{f:digraph}.
  \tikzset{middlearrow/.style={
        decoration={markings,
            mark= at position 0.5 with {\arrow{#1}} ,
        },
        postaction={decorate}
    }
}
\begin{figure}[h]
\begin{center}
\begin{tikzpicture}
\node[] (1) at (0,0) {$v_1$};
    \node[] (2) at (1.5,0) {$v_2$};
    \node[] (3) at (0.7,-1) {$v_3$};
    \node[] (4) at (2.5,-1) {$v_4$};
    \node[] (5) at (3.3,-2.0) {$v_5$};
    \node[] (6) at (4,-1) {$v_6$};
    \node[] (7) at (6, 0) {$v_7$};
    \node[] (8) at (5,0) {$v_8$};
    \node[] (9) at (3.3,-3) {$v_9$};
    \node[] (10) at (3.3,-4) {$v_{10}$};
     \node[draw=black,fill=pink] (11) at (-1,-2) {$u_1$};
     \node[draw=black,fill=pink] (12) at (6,-2) {$u_2$};
     \node[draw=black,fill=pink] (13) at (4.7,-3) {$u_3$};
 
     \path [->,solid, cyan, line width=1.0pt](11) edge node[right] {} (1); 
     \path [->,solid, cyan, line width=1.0pt](11) edge node[right] {} (3);
     \path [->,solid, cyan, line width=1.0pt](11) edge node[right] {} (4);
     \path [->,solid, cyan, line width=1.0pt](12) edge node[right] {} (6);
     \path [->,solid, cyan, line width=1.0pt](12) edge node[right] {} (7);
     \path [->,solid, cyan, line width=1.0pt](13) edge node[right] {} (8);
     \path [->,solid, cyan, line width=1.0pt](13) edge node[right] {} (10);

     \draw[bend right,middlearrow={>}, line width=1.0pt]   (2) to node [auto] {} (1);
     \draw[bend right,middlearrow={>}, line width=1.0pt]   (1) to node [auto] {} (3);
     \draw[bend right,middlearrow={>}, line width=1.0pt]   (3) to node [auto] {} (2);
     \draw[bend right,middlearrow={>}, line width=1.0pt]   (9) to node [auto] {} (10);
     \draw[bend right,middlearrow={>}, line width=1.0pt]   (10) to node [auto] {} (9);
     \draw[bend right,middlearrow={>}, line width=1.0pt]   (7) to node [auto] {} (8);
     \draw[bend right,middlearrow={>}, line width=1.0pt]   (8) to node [auto] {} (7);
     \draw[bend right,middlearrow={>}, line width=1.0pt]   (6) to node [auto] {} (4);
     \draw[bend right,middlearrow={>}, line width=1.0pt]   (5) to node [auto] {} (6);
     \draw[bend right,middlearrow={>}, line width=1.0pt]   (4) to node [auto] {} (5);
     \path [->,solid, line width=1.0pt](3) edge node[right] {} (4);
     \path [->,solid, line width=1.0pt](8) edge node[right] {} (6);
     \path [->,solid, line width=1.0pt](9) edge node[right] {} (5);
     \path(4) edge [loop above, line width=1.0pt] node {} (4);
      \path(2) edge [loop above, line width=1.0pt] node {} (2);   
   \end{tikzpicture}
\end{center}
\caption{Illustration of the digraph \(G(A,B)\). Each vertex in a coloured box represents an input. The black edges denote the state-connections and the cyan coloured edges denote input-connections.}
\label{f:digraph}
\end{figure}
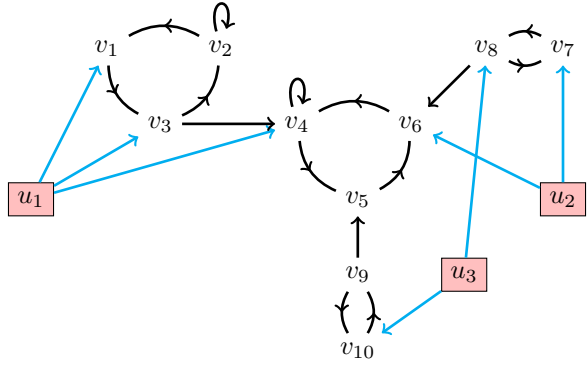

Observe that \((A, B)\) is structurally controllable. Recall from \S\ref{ss:problems} that \(w_{ij}\) is the cost of the input-connection between the input \(u_j\) and the state \(v_i\). Let the costs in the present case be \(w_{11}=15\), \(w_{31}=10\), \(w_{41}=20\), \(w_{62}=15\), \(w_{72}=5\), \(w_{83}=5\), and \(w_{10,3}=10\). Clearly, \(\BA\) has a perfect matching and the SSCCs of \(G(A)\) are: \(\mathcal{S}_1=\lbrace v_1, v_2, v_3 \rbrace\), \(\mathcal{S}_2=\lbrace v_7, v_8 \rbrace\), and \(\mathcal{S}_3=\lbrace v_9, v_{10} \rbrace\). 

We resort to Algorithm \ref{AlgorithmProblem2} to solve \eqref{P2}. After we execute Algorithm \ref{AlgorithmProblem2}, the solution obtained is \(B^{*}\) with \(B^{*}_{31}=\star\), \(B^{*}_{72}=\star\), and \(B^{*}_{10,3}=\star\), and it has the (minimum) cost of \(25\). Clearly, \(G(A, B^{*})\) has all the SSCCs \(\mathcal{S}_1\), \(\mathcal{S}_2\), and \(\mathcal{S}_3\) of \(G(A)\) accessible from some input vertex, and the sum of the cost of the input-connections in \(B^{*}\) have the least value among all the input matrices \(B^{\prime}\in \mathcal{K}\). It follows from Proposition \ref{P:problem3assum1} that the input matrix \(B^{*}\) calculated as above also solves \eqref{P3}.

If uniform non-zero cost is assumed for all input-connections then we get an input matrix \(B^{*}\) with \(B^{*}_{11}=\star\), \(B^{*}_{83}=\star\), and \(B^{*}_{10,3}=\star\) that solves \eqref{P1}. Observe that all the SSCCs \(\mathcal{S}_1\), \(\mathcal{S}_2\), and \(\mathcal{S}_3\) of \(G(A)\) are accessible from the input vertices in \(G(A, B^{*})\); of course, the input vertices associated with each SSCCs may not be distinct from each other.
\\

\noindent \textit{Example 2}: Let\\
\begin{equation}
\label{example2}
A=\begin{bmatrix}
0 & \star & 0 &  0 & 0 & 0 & 0 & 0  \\
0 & 0 & \star  &  0 & 0 & 0 & \star  & 0 \\
0 & 0 & 0 &  \star  & 0 & 0 & 0 & 0 \\
\star  & 0 & 0 &  0 & 0 & 0 & 0 & \star  \\
\star  & 0 & 0 &  0 & 0 & \star  & 0 & 0\\
0 & 0 & 0 &  0 & \star  & 0 & 0 & 0\\
0 & 0 & 0 &  0 & \star  & 0 & 0 & 0\\
0 & 0 & 0 &  0 & \star  & 0 & 0 & 0\\
\end{bmatrix}
,
B=\begin{bmatrix}
0 & 0  & 0 & \star \\
\star  & 0  & 0 & 0\\
0  & \star   & 0  & 0 \\
0 & 0  & 0 & 0\\
0 & 0  & \star & 0\\
\star  & \star   & \star  & \star \\
\star  & \star   & \star  & \star \\
\star  & \star   & \star  & \star \\
\end{bmatrix}.
\end{equation}
\begin{figure}[th]
\begin{center}
\begin{tikzpicture}
\node[] (1) at (4,0) {$v_1^1$};
    \node[] (2) at (4,-0.5) {$v_2^1$};
    \node[] (3) at (4,-1) {$v_3^1$};
    \node[] (4) at (4,-1.5) {$v_4^1$};
     \node[] (5) at (4,-2) {$v_5^1$};
     \node[] (6) at (4,-2.5) {$v_6^1$};
     \node[] (7) at (4,-3) {$v_7^1$};
     \node[] (8) at (4,-3.5) {$v_8^1$};
    \node[] (1') at (7.5,0) {$v_1^2$};
    \node[] (2') at (7.5,-0.5) {$v_2^2$};
    \node[] (3') at (7.5,-1) {$v_3^2$};
    \node[] (4') at (7.5,-1.5) {$v_4^2$};
    \node[] (5') at (7.5,-2) {$v_5^2$};
     \node[] (6') at (7.5,-2.5) {$v_6^2$};
     \node[] (7') at (7.5,-3) {$v_7^2$};
     \node[] (8') at (7.5,-3.5) {$v_8^2$};
     \node[] (u1) at (4,-4) {$u_1$};
     \node[] (u2) at (4,-4.5) {$u_2$};
     \node[] (u3) at (4,-5) {$u_3$};
     \node[] (u4) at (4,-5.5) {$u_4$};
\draw (2) edge[-,solid,draw=black, line width=1.0pt] (1') ;
\draw (3) edge[-,solid,draw=black, line width=1.0pt] (2') ;
\draw (7) edge[-,solid,draw=black, line width=1.0pt] (2') ;
\draw (4) edge[-,solid,draw=black, line width=1.0pt] (3') ;
\draw (1) edge[-,solid,draw=black, line width=1.0pt] (4') ;
\draw (8) edge[-,solid,draw=black, line width=1.0pt] (4') ;
\draw (1) edge[-,solid,draw=black, line width=1.0pt] (5') ;
\draw (6) edge[-,solid,draw=black, line width=1.0pt] (5') ;
 \draw (5) edge[-,solid,draw=black, line width=1.0pt] (6') ;
\draw (5) edge[-,solid,draw=black, line width=1.0pt] (7') ;
\draw (5) edge[-,solid,draw=black, line width=1.0pt] (8') ;
\draw (u1) edge[-,draw=cyan, line width=1.0pt] (2') ;
\draw (u1) edge[-,draw=cyan, line width=1.0pt] (6') ;
\draw (u1) edge[-, draw=cyan, line width=1.0pt] (7') ;
\draw (u1) edge[-,draw=cyan, line width=1.0pt] (8') ;
\draw (u2) edge[-,draw=cyan, line width=1.0pt] (3') ;
\draw (u2) edge[-,draw=cyan, line width=1.0pt] (6') ;
\draw (u2) edge[-,draw=cyan, line width=1.0pt] (7') ;
\draw (u2) edge[-,draw=cyan, line width=1.0pt] (8') ;
\draw (u3) edge[-,draw=cyan, line width=1.0pt] (5') ;
\draw (u3) edge[-,draw=cyan, line width=1.0pt] (6') ;
\draw (u3) edge[-,draw=cyan, line width=1.0pt] (7') ;
\draw (u3) edge[-,draw=cyan, line width=1.0pt] (8') ;
\draw (u4) edge[-,draw=cyan, line width=1.0pt] (1') ;
\draw (u4) edge[-,draw=cyan, line width=1.0pt] (6') ;
\draw (u4) edge[-,draw=cyan, line width=1.0pt] (7') ;
\draw (u4) edge[-,draw=cyan, line width=1.0pt] (8') ;
\end{tikzpicture}
\end{center}
\caption{Illustration of the bipartite graph \(\BAB\) associated with the pair \((A, B)\). The black edges corresponds to the state-connections. The cyan coloured edges corresponds to the input-connection.}
\label{f:bipartitegraph}
\end{figure}
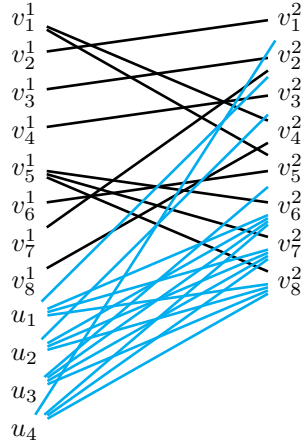

\noindent Observe that \((A, B)\) is structurally controllable. Suppose that costs of the input-connections are \(w_{21}=5\), \(w_{61}=1\), \(w_{71}=2\), \(w_{81}=2\), \(w_{32}=5\), \(w_{62}=1\), \(w_{72}=5\), \(w_{82}=6\), \(w_{53}=5\), \(w_{63}=2\), \(w_{73}=1\), \(w_{83}=2\), \(w_{14}=5\), \(w_{64}=5\), \(w_{74}=3\), and \(w_{84}=1\). \(G(A)\) is strongly connected, and therefore we resort to Algorithm \ref{algorithmstronglyprob2} to solve \eqref{P3}. The bipartite graph \(\BAB\) is depicted in Fig.~\ref{f:bipartitegraph}. A MCMM \(M_1\) of total cost \(2\) in \((\BAB;c)\) (under the cost function \(c\) defined in \eqref{e:weightProblem3assum4.2}) is \(\lbrace (v_1^1, v_4^2)\),\((v_4^1,v_3^2),(v_3^1, v_2^2),(v_2^1, v_1^2),(v_5^1,v_6^2) (v_6^1, v_5^2) ,\) 
\((u_3, v_7^2),(u_4,v_8^2)\rbrace\). Therefore, the solution we obtain for \eqref{P3} is \(B^{*}\) with \(B^{*}_{73}=\star\), and \(B^{*}_{84}=\star\). Since each input-connection in \(G(A, B)\) has positive cost, \(B^{*}\) also solves \eqref{P2}.

If all the costs are uniform, then we may obtain another MCMM \(M_2\) in \((\BAB;c)\) (under the cost function \(c\) defined in \eqref{e:weightProblem1assum4.2}) is \(M_2=\lbrace (v_1^1, v_5^2),(v_2^1,v_1^2),(v_3^1, v_2^2),(v_4^1, v_3^2),\)
\((v_8^1,v_4^2) (v_5^1, v_8^2),(u_1, v_6^2),\) \((u_2,v_7^2)\rbrace\) and \(B^{*}\) has \(\star\) entries at \(B^{*}_{61}=\star\), and \(B^{*}_{72}=\star\), which is a solution of \eqref{P1}.

\noindent \textit{Example 3}:
We provide an illustration of our results in the context of the benchmark IEEE 118-bus system. This particular dynamical system corresponds to an electric power grid composed of
\begin{itemize}[label=\(\circ\), leftmargin=*]
\item  \(118\) buses;
\item \(53\) power generators, and,
\item \(65\) power loads,
\end{itemize} 
all coupled via power lines. We adopt the cyber-physical model of the generators and loads  proposed in \cite{ref:IliXieKhaMou-10}, where a linear system was obtained via Taylor linearization at the nominal operating point. This linear model yields the following: the total number of vertices in \(G(A)\) is \(407\) and the total number of directed edges is \(920\). 

The dynamics of the electric power grid and its components, modeled as in \cite{ref:IliXieKhaMou-10}, is described by the following
state variables: \(P_{T_{G_i}}\) represents the mechanical power of the
turbine of the generator \(G_i\), \(w_{G_i}\) the generator \(G_i\) frequency and \(a_{G_i}\) its valve opening. In addition, \(I_{L_j}\) is the real energy consumed by the load \(L_j\) and \(w_{L_j}\) the frequency measured at load \(L_j\) location. The different components are connected
through the injected/received power to/from the network at the connection site, which dynamics depend on the frequency of the components on the neighbouring buses; the injected and received power variables for generator \(i\) and load \(j\) are \(P_{G_i}\) and \(P_{L_j}\), respectively. If the generator and load are annexed to a bus, hence to the network, then the induced dynamics is coupled with the power injected/received to/from
the network. Thus, the digraph representation of the dynamics has bidirectional connection between the injected/received power variables and the frequency of the corresponding components (corresponding to the injected/received power to/from the network), as depicted in Fig.~\ref{f:genload}.

\tikzset{middlearrow/.style={
        decoration={markings,
            mark= at position 0.5 with {\arrow{#1}} ,
        },
        postaction={decorate}
    }
} 

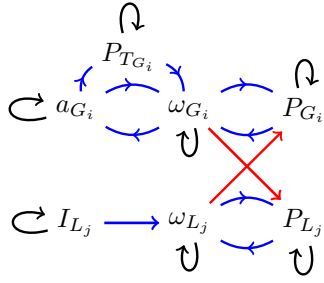
\begin{figure}[h!]
\begin{center}
\begin{tikzpicture}
\node[] (0) at (0,0) {\(a_{G_i}\)};
\node[] (1) at (1.5, 0) {\(\omega_{G_i}\)};
 \node[] (2) at (3,0) {\(P_{G_i}\)};
 \node[] (3) at (3,-1.5) {\(P_{L_j}\)};
 \node[] (4) at (1.5,-1.5) {\(\omega_{L_j}\)};
 \node[] (5) at (0,-1.5) {\(I_{L_j}\)};
  \node[] (6) at (0.7,0.7) {\(P_{T_{G_i}}\)};
  \draw (1) edge[->,draw=red, line width=1.0pt] (3) (4) edge[->,draw=red, line width=0.9pt] (2) (5) edge[->,draw=blue, line width=1.0pt] (4); 
      \draw[bend left,middlearrow={>},draw=blue, line width=0.9pt] (0) to node [auto] {} (1);
   \draw[bend left,middlearrow={>},draw=blue, line width=0.9pt] (1) to node [auto] {} (0);
   \draw[bend left,middlearrow={>},draw=blue, line width=0.9pt] (1) to node [auto] {} (2);
  \draw[bend left,middlearrow={>},draw=blue, line width=0.9pt] (2) to node [auto] {} (1);
  \draw[bend left,middlearrow={>},draw=blue, line width=0.9pt] (3) to node [auto] {} (4);
  
   \draw[bend left,middlearrow={>},draw=blue, line width=0.9pt] (4) to node [auto] {} (3);
   
   \draw[bend left,middlearrow={>},draw=blue, line width=0.9pt] (0) to node [auto] {}
    (6);
    \draw[bend left,middlearrow={>},draw=blue, line width=0.9pt] (6) to node [auto] {} (1);
  \path(1) edge [loop below, line width=0.9pt] node {} (1);
    \path(2) edge [loop above, line width=0.9pt] node {} (2);
  \path(0) edge [loop left, line width=0.9pt] node {} (0);
  \path(4) edge [loop below, line width=0.9pt] node {} (4);
  \path(3) edge [loop below, line width=0.9pt] node {} (3);
  \path(5) edge [loop left, line width=0.9pt] node {} (5);
   \path(6) edge [loop above, line width=0.9pt] node {} (6);

    \end{tikzpicture}
   \end{center}
   \caption{Illustration of coupling between two neighbouring components, a generator \(i\) connected to a load \(j\) through transmission line \((i, j)\).} 
 \label{f:genload}
\end{figure} 

Let us assume that generator \(i\) and load \(j\) are attached to the same bus, or different buses but there exists a transmission line \((i,j)\) between them. The frequency component of bus \(i\), i.e., \(w_{G_i}\) affects the dynamics of the power component of bus \(j\), i.e., \(P_{L_j}\), and vice-versa. This implies that there exist outgoing  edges  from  the  frequencies  of  the components into the powers of the components in the neighbouring  buses, as shown in Fig.~\ref{f:genload}.

\begin{figure}[h!]
\centering
\includegraphics[width = 3.2 in, keepaspectratio]{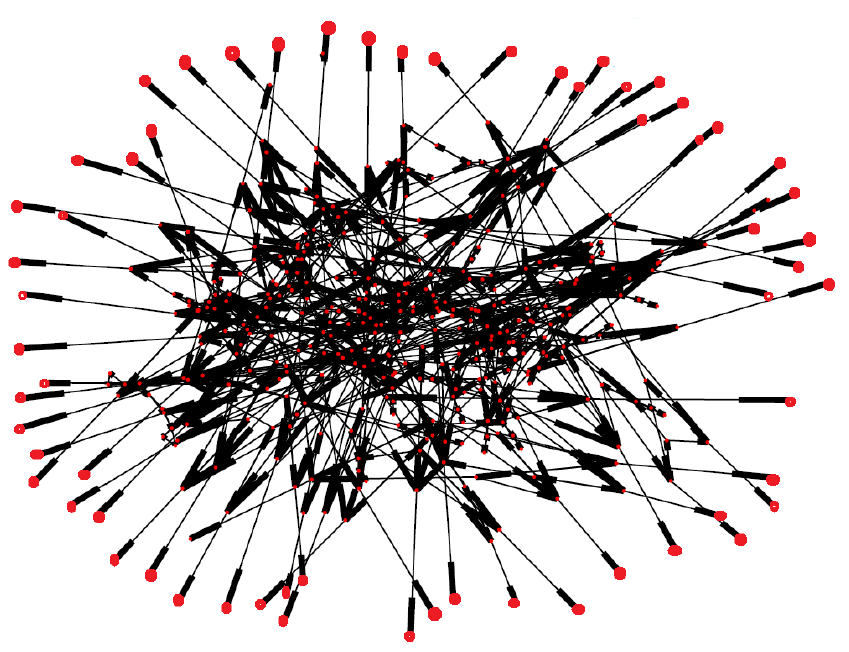}
\caption{Graphical representation of the IEEE 118-bus system, \(G(A)\). The SSCCs of \(G(A)\) are depicted in bold (red).}
\label{f:plot118}
\end{figure}

We provide a solution of \eqref{P2} for the IEEE 118-bus system. Note that by construction, as shown in Fig. \ref{f:genload}, each vertex of the digraph has a self loop. Therefore, Assumption \ref{a:perfectmatching} holds for \(\BA\). It is proved in \cite[Theorem 2]{ref:RamPeqAguKar-15} that if the network topology is connected then the number of SSCCs present in \(G(A)\) associated with a power grid is equal to the number of loads. Clearly, if there exists a load \(j\), then there is an SSCC given by \(I_{L_j}\). Since the topology of IEEE \(118\)-bus system is connected, the number of SSCCs is \(65\), i.e., equal to the number of loads. The graphical representation of IEEE 118-bus system with the SSCCs represented in bold (red) is shown in Fig.~\ref{f:plot118}.  For simplicity the self loops are omitted.

\begin{figure}[h!]
\centering
\includegraphics[width = 3.5 in, keepaspectratio]{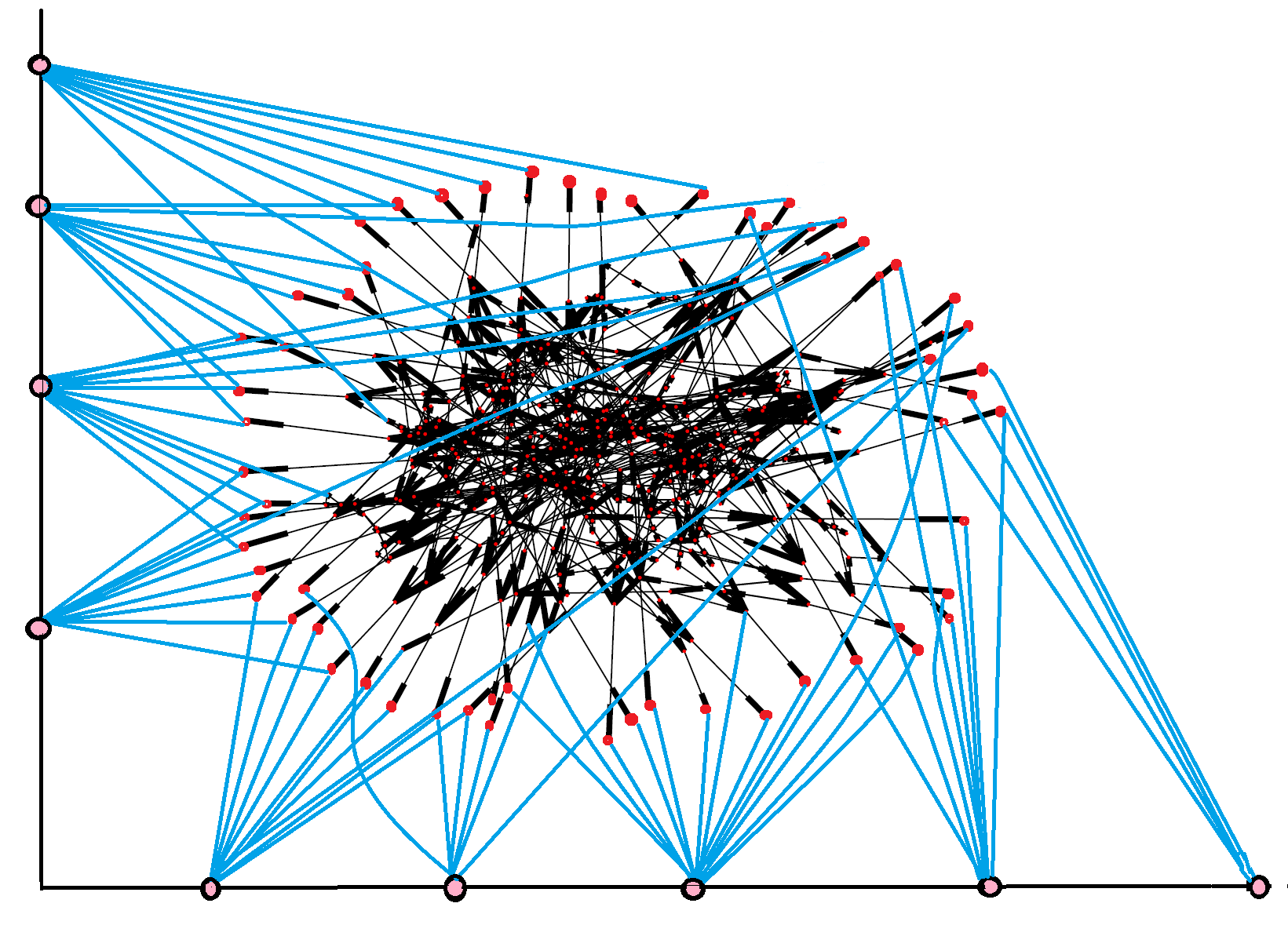}
\caption{Illustration of the digraph associated with IEEE 118-bus system along with the set of inputs and input-connections (coloured in cyan)  are depicted in the figure. The black-coloured edges represent state-connections. The pink coloured vertices with black border represents the inputs.}
\label{f:IEEE118input}
\end{figure}

For IEEE 118-bus system, the given set of inputs and the input-connections associated are depicted in Fig.~\ref{f:IEEE118input} such that the given system is structurally controllable.  All the input-connections are assumed to have cost in the range \((0,20)\). For simplicity, it is not depicted in the figure.

The following Fig.~\ref{f:IEEEsolution} depicts a solution of \eqref{P2} by using Algorithm \ref{AlgorithmProblem2}.  Observe that Fig.~\ref{f:IEEEsolution} has fewer input-connections that Fig.~\ref{f:IEEE118input}. We obtain a set of input-connections of the least cost needed to ensure structural controllability, i.e., each \(I_{L_j}\) has an input-connection from some input in Fig.~\ref{f:IEEEsolution}. Notice that from a physical point of view the actuation of the variables \(I_{L_j}\) in each load means that we need to be able to actuate the real power consumed by the aggregate load. 
\begin{figure}[h!]
\centering
\includegraphics[width = 3.5 in, keepaspectratio]{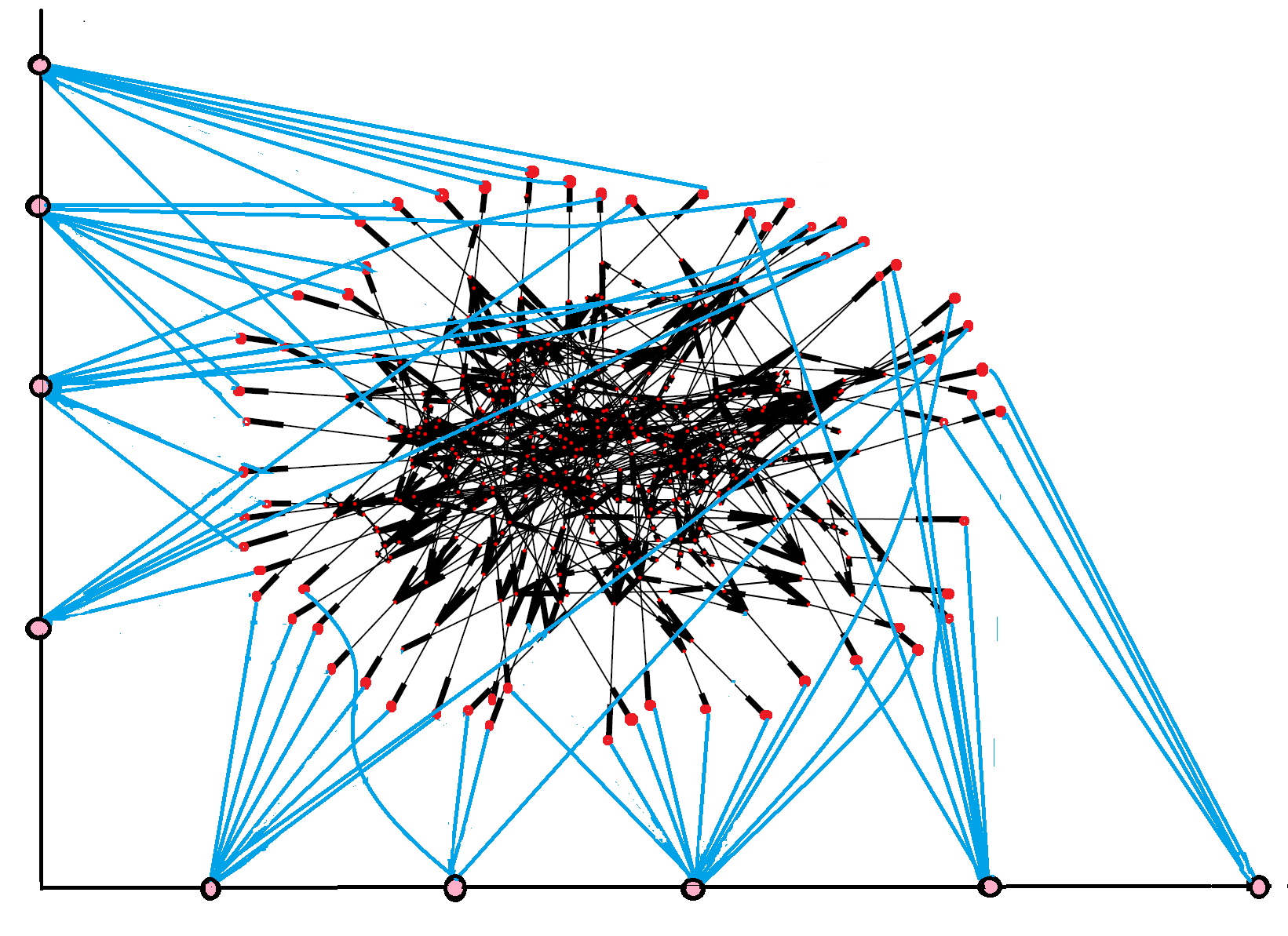}
\caption{Illustration of the digraph associated with IEEE 118-bus system along with the set of inputs and input-connections (coloured in cyan) retained in the solution of \eqref{P2} are depicted in the figure. The black-coloured edges represent state-connections. The pink coloured vertices with black border represents the inputs.}
\label{f:IEEEsolution}
\end{figure}

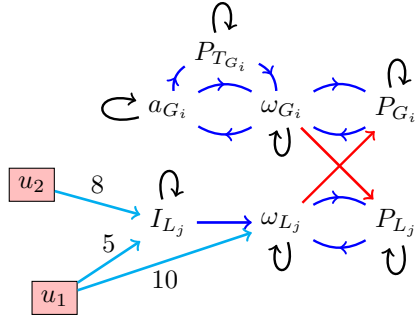
\begin{figure}[h!]
\begin{center}
\begin{tikzpicture}
\node[] (0) at (0,0) {\(a_{G_i}\)};
\node[] (1) at (1.5, 0) {\(\omega_{G_i}\)};
 \node[] (2) at (3,0) {\(P_{G_i}\)};
 \node[] (3) at (3,-1.5) {\(P_{L_j}\)};
 \node[] (4) at (1.5,-1.5) {\(\omega_{L_j}\)};
 \node[] (5) at (0,-1.5) {\(I_{L_j}\)};
  \node[] (6) at (0.7,0.7) {\(P_{T_{G_i}}\)};
   \node[draw=black,fill=pink] (7) at (-1.5,-2.5) {$u_1$};
     \node[draw=black,fill=pink] (8) at (-1.8,-1) {$u_2$};
   \draw (1) edge[->,draw=red, line width=1.0pt] (3) (4) edge[->,draw=red, line width=0.9pt] (2) (5) edge[->,draw=blue, line width=1.0pt] (4); 
      \draw[bend left,middlearrow={>},draw=blue, line width=0.9pt] (0) to node [auto] {} (1);
   \draw[bend left,middlearrow={>},draw=blue, line width=0.9pt] (1) to node [auto] {} (0);
   \draw[bend left,middlearrow={>},draw=blue, line width=0.9pt] (1) to node [auto] {} (2);
  \draw[bend left,middlearrow={>},draw=blue, line width=0.9pt] (2) to node [auto] {} (1);
  \draw[bend left,middlearrow={>},draw=blue, line width=0.9pt] (3) to node [auto] {} (4);
  
   \draw[bend left,middlearrow={>},draw=blue, line width=0.9pt] (4) to node [auto] {} (3);
   
   \draw[bend left,middlearrow={>},draw=blue, line width=0.9pt] (0) to node [auto] {}
    (6);
    \draw[bend left,middlearrow={>},draw=blue, line width=0.9pt] (6) to node [auto] {} (1);
  \path(1) edge [loop below, line width=0.9pt] node {} (1);
    \path(2) edge [loop above, line width=0.9pt] node {} (2);
  \path(0) edge [loop left, line width=0.9pt] node {} (0);
  \path(4) edge [loop below, line width=0.9pt] node {} (4);
  \path(3) edge [loop below, line width=0.9pt] node {} (3);
  \path(5) edge [loop above, line width=0.9pt] node {} (5);
   \path(6) edge [loop above, line width=0.9pt] node {} (6);
   \path [->,draw=cyan, line width=1.0pt](7) edge node[right,text width=1cm,midway,above,align=center] {5} (5); 
     \path [->,draw=cyan, line width=1.0pt](7) edge node[right,text width=1cm,midway,below,align=center] {10} (4); 
     \path [->,draw=cyan, line width=1.0pt](8) edge node[right,text width=1cm,midway,above,align=center] {8} (5);      
    \end{tikzpicture}
   \end{center}
   \caption{Illustration of digraph representation of generator \(i\) connected to load \(j\) through transmission line \((i, j)\) and the input-connections (coloured in cyan) associated with given inputs. Each input-connection has a cost depicted by the number over it.}
 \label{f:genloadinputs}
\end{figure}

The number of states in the linearized model of IEEE-118 bus network is roughly 400. For the ease of understanding, we depict the use of Algorithm \ref{AlgorithmProblem2} by zooming in on a very small part of the  network which consists of a generator and load connected by a transmission line, as shown in Fig. \ref{f:genload}. The generator \(i\) and the load \(j\) along with the set of input-connections (coloured in cyan) are depicted in Fig~\ref{f:genloadinputs}. Each input-connection has a cost associated shown by the number over it in the figure. 
 
The SSCC in the digraph shown in Fig.~\ref{f:genloadinputs} is \(I_{{L}_j}\). Note that \(I_{{L}_j}\) form a SSCC of the digraph \(G(A)\) also. After we execute the Algorithm \ref{AlgorithmProblem2} only the input-connection connecting input \(u_1\) to state \(I_{{L}_j}\), i.e., \((u_1, I_{{L}_j})\) is retained with total cost of \(5\). 

\section{Appendix}
\label{s:appendix}
\textit{Proof of Lemma \ref{l:sparseBperfectmatching}}: We establish the assertion in two steps: In step (i) we prove that \(\norm{B^*}_0 \leq q\) and in step (ii) we show that \(\norm{B^*}_0 \geq q\).

Step (i). Since Assumption \ref{a:perfectmatching} holds, the pair \((A,B)\) is structurally controllable if and only if all the SSCCs of \(G(A)\) are accessible from the input vertices. Therefore, every SSCCs has at least one state vertex directly connected to one of the input vertices in \(G(A,B)\). Also, for each  SSCC exactly one input-connection is sufficient to ensure accessibility. This confirms that there exists a \(B^{\prime}\) such that \((A,B^{\prime})\) is structurally controllable, i.e., \(B^{\prime}\in \mathcal{K}\) and  \(\norm{B^{\prime}}_0=q\). In other words, \(\norm{B^{*}}_0 \leq q\).

Step (ii). Suppose that there exists a \(B^{\prime} \in \mathcal{K}\) such that \(\norm{B^{\prime}}_0 < q\). We assume that \(\norm{B^{\prime}}_0 =q-1\). It means that there are \(q\) SSCCs in \(G(A)\) and only \(q-1\) input-connections. Since all the SSCCs are vertex-disjoint from each other there exists at least one SSCC not accessible from any input vertex. This contradicts the assumption that \((A, B^{\prime})\) is structurally controllable. Therefore, every \(B^{\prime\prime} \in \mathcal{K}\) is such that  \(\norm{B^{\prime\prime}}_0 \geq q\), leading to \(\norm{B^{*}}_0 \geq q\). The assertion follows.
\qed
\\

\textit{Proof of Theorem \ref{t:problemP2}}: It follows from the procedure outlined in Algorithm \ref{AlgorithmProblem2} that exactly one state vertex, having an input-connection from some input vertex, is selected for each SSCC at each iteration. Therefore, \((A, B^{*})\) is structurally controllable, i.e., \(B^{*} \in \mathcal{K}\) and \(\norm{B^{*}}_0=q\). By Lemma \ref{l:sparseBperfectmatching} it follows that the \(B^{*}\) so obtained has minimum number of non-zero entries. The procedure employed in Algorithm \ref{AlgorithmProblem2} also ensures that \(B^{*}\) has the least cost among the collection of all sparsest \(B^{\prime}\in \mathcal{K}\), i.e., \(\norm{B^{*}}_w= \sum_{\ell=1}^{d}\sum_{k=1}^{m}w_{\ell k} \textbf{1}_{\lbrace B_{\ell k}^{*}\neq 0 \rbrace}\) has the least value. 
\qed
\\

\textit{Proof of Proposition \ref{P:problem3assum1}}: Suppose that the assertion is false, and there exists another \(B^{\prime}\in \mathcal{K}\) such that \( \norm{B^{\prime}}_w < \norm{B^{*}}_w\). If \(\norm{B^{\prime}}_0 = q =\norm{B^{*}}_0\) then our assumption is false. So, consider the case, where \(\norm{B^{\prime}}_0 > q =\norm{B^{*}}_0\). Without loss of generality, assume that \(\norm{B^{\prime}}_0 = q+1\). Since \(B^{\prime} \in \mathcal{K}\) implies that at least \(q\) of the input-connections are connected to one state vertex in each SSCC. Therefore, it is possible to extract a new input matrix \(B^{\prime\prime} \in \mathcal{K}\)  from \(B^{\prime}\) such that \(\norm{B^{\prime\prime} }_0 = q\) and \(\norm{B^{\prime\prime} }_w \leq \norm{B^{\prime}}_w < \norm{B^{*}}_w\). This contradicts the optimality of \(B^{*}\), and completes the proof.  
\qed
\\

\textit{Proof of Theorem \ref{t:assumption4.2Prob2}}: It is given that pair \((A, B)\) is structurally controllable. Proposition \ref{p:Aols}  implies that there exists an SDR for \(\bir\) in \(\BAB\). The cost structure (see Step 3) of Algorithm \ref{algorithmstronglyprob2} ensures that when a MCMM is computed, it uses the input-connections which minimizes the cost and obtains an SDR for \(\bir\) in \(\BAB\). Therefore, the obtained input matrix \(B^{*}\) (see Step 5 of Algorithm \ref{algorithmstronglyprob2}) is a solution of \eqref{P3}.

Given the pair \((A, B)\), the construction of \(\BAB\) has linear complexity. The problem of finding a MCMM in \(\BAB\) can be efficiently solved using Hungarian algorithm \cite{ref:JMun-57} with computation complexity of \(O\big((d + m)^3\big)\) under the cost function \(c\) defined in \eqref{e:weightProblem3assum4.2}, where \(d\) and \(m\) are the number of state and input vertices in \(G(A,B)\). Therefore, the complexity of Algorithm \ref{algorithmstronglyprob2} is dominated by Step 4 implying an overall complexity of \(O((d+m)^3)\), as asserted.
\qed
\\

\textit{Proof of Lemma \ref{l:solveP2P3}}: It is well-known that a matching in a graph is maximum if and only if the graph has no augmenting path \cite[Theorem 3.1.10, p.~109]{ref:DouWesT-96}. \footnote{Given a matching M in \(\BA\), an M-augmenting path is a path that alternates between edges in M and not in M and whose endpoints are unmatched by M.} Let \(M^{*}\) be the MCMM (obtained in Algorithm \ref{algorithmstronglyprob2}) associated with \(B^*\) saturating the vertices of \(\bir\). Let \(M^{\prime}= M^{*} \cap \mathcal{E}_A\). Note that \(B^{*}\) solves \eqref{P3}. To show that \(B^{*}\) also solves \eqref{P2} it suffices to show that \(B^{*}\) has the least number of non-zero entries. We prove this by showing that \(M^{\prime}\) is a maximum matching in \(\BA\). Note that if \(M^{\prime}\) is a maximum matching in \(\BA\) then the additional number of edges of \(\mathcal{E}_B\) required to obtain a matching saturating all the vertices of \(\bir\) in \(\BAB\) is necessarily the least among all saturating matchings in \(\BAB\). 

Next, we prove that \(M^{\prime}\) is a maximum matching in \(\BA\). Suppose that the assertion is false, then there exists an \(M^{\prime}\)-augmenting path in \(\BA\). We use this \(M^{\prime}\)-augmenting path to obtain a new matching \(M^{\prime\prime}\) in \(\BA\) whose cardinality is exactly one greater than \(M^{\prime}\). Clearly, \(M^{\prime\prime}\) saturates hitherto an unmatched vertex \(v_r^2\) by \(M^{\prime}\) of \(\bir\). Since \(M^{*}\) saturates all the vertices of \(\bir\) in \(\BAB\) the vertex \(v_r^2\) is saturated by some input vertex, say \(u_j \in V_B\), i.e., \((u_j, v_r^2)\in M^{*}\). Consider \(M_1=M^{\prime\prime}\sqcup \big\lbrace (M^{*}\cap \mathcal{E}_B)\setminus (u_j,v_r^2) \big\rbrace\). Notice that \(M_1\) is a matching in \(\BAB\) that saturates all the vertices of \(\bir\). Thus, \(\abs{M^{*}}=\abs{M_1}=d\). Also, the cost of matching \(M_1\) is \(\sum_{e \in M_1} c(e)= \sum_{e \in M^{*}} c(e) -w_{rj}+0\). Since \(w_{rj}> 0\), cost of \(M_1\) is strictly less than cost of \(M^{*}\). This contradicts that \(M^{*}\) is a MCMM in \(\BAB\) and completes the proof. Thus, the obtained \(B^{*} \in \mathcal{K}\) is minimal.    
\qed
\\

\textit{Proof of Lemma \ref{o:c1andc2}}:
Since \(\BAB\) admits an SDR any maximum matching of \(\BAB\) saturates all the vertices of \(\bir\). For a matching \(M\), let \(c(M)= \sum_{e \in M} c(e)\) denote the cost of \(M\) under cost function \(c\).

Let \(M^{*}\) be a MCMM under \(c_1\) in \(\BAB\). We know that \(M^{*}=(M^{*}\cap \mathcal{E}_A)\sqcup (M^{*}\cap \mathcal{E}_B)\). Since the edges of \(\mathcal{E}_B\) have positive cost under \(c_1\), \(M^{*}\cap \mathcal{E}_A\) is a maximum matching in \(\BA\)  by Lemma \ref{l:solveP2P3}. Let \(\ell\) be the number of vertices of \(\bir\) unmatched by \(M^{*}\cap \mathcal{E}_A\). Then \(\abs{M^{*}\cap \mathcal{E}_B}=\abs{M^{*}}-\abs{M^{*}\cap \mathcal{E}_A}=\ell\) and \(c(M^{*})=c_1(M^{*})-\ell\). To show that \(M^{*}\) is a MCMM under cost function \(c\) we proceed by contradiction. Suppose there exist another maximum matching \(M\) in \(\BAB\) such that \(c(M)<c(M^{*})\). 

Now, \(M=(M\cap \mathcal{E}_A)\sqcup (M\cap \mathcal{E}_B)\) and \(M\cap \mathcal{E}_A\) is a (possibly not maximum) matching in \(\BA\). Then \(\abs{M\cap \mathcal{E}_B}\geq \ell\). If \(S\) be the set of vertices of \(\BA\) saturated by \(M\cap \mathcal{E}_A\) then there exists some maximum matching \(M^{\prime}\) in \(\BA\) saturating all of \(S\) \cite[Chapter 3, p.~118]{ref:DouWesT-96}. Let \(\mathcal{E}^{\prime}=\big\lbrace (u_j, v_i^2)\in M \;|\; \exists\; v_k^1\in \bil\; \text{such that}\; (v_k^1, v_i^2) \in M^{\prime}\big\rbrace\). Let \(M^{\prime\prime}=M^{\prime}\sqcup \big \lbrace (M\cap \mathcal{E}_B) \setminus \mathcal{E}^{\prime}\big \rbrace\). Note that \(\abs{M^{\prime\prime}}=\abs{M}\) and \(c(M^{\prime\prime})\leq c(M)\). Then,
\begin{equation*}
\begin{aligned}
 c_1(M^{\prime\prime})&=c(M^{\prime\prime})+\ell\leq c(M)+\ell &<c(M^{*})+\ell & = c_1(M^{*}).
\end{aligned}
\end{equation*}
This contradicts that \(M^{*}\) is a MCMM in \(\BAB\) under cost function \(c_1\) and completes the proof.
\qed
\\

\textit{Proof of Proposition \ref{p:directedtreeP3}}: Let \(M^{*}\) be the MCMM associated with \(B^{*}\). Since \(G(A)\) has exactly one SSCC \(\mathcal{S}_1=\lbrace v_i\rbrace\) and \(\deg^{-}{v_i}=0\) in \(G(A)\), \(M^{*}\) has an input-connection of the form \((u_j, v_i^2)\in M^{*}\cap \mathcal{E}_B\) which satisfy the accessibility condition needed for structural controllability as well as saturates \(v_i\) in \(\Gamma(A, B^{*})\). The rest of the input-connections corresponding to \(M^{*}\cap \mathcal{E}_B\) are required to satisfy only the no-dilation criterion. Thus, \(B^{*}\in \mathcal{K}\). We prove that \(B^{*}\) has the least cost by contradiction. Suppose there exists another \(B^{\prime} \in \mathcal{K}\) such that \(\norm{B^{\prime}}_w < \norm{B^{*}}_w\). Let \(M^{\prime}\) be a matching saturating \(\bir\) in \(\Gamma(A, B^{\prime})\). Note that \(M^{\prime}\) must have an input-connection of the form \((u_j, v_i^2)\) in \(\Gamma(A, B^{\prime})\). Then \(\sum_{e\in M^{\prime}}c(e)\leq \norm{B^{\prime}}_w < \norm{B^{*}}_w=\sum_{e\in M^{*}}c(e)\), contradicting the optimality of \(M^{*}\) and completing the proof.
\qed
\\

\textit{Proof of Theorem \ref{t:approximate}}: The obtained \(B^{*} \in \mathcal{K}\), since it ensures that all the SSCCs are accessible in \(G(A, B^{*})\) and an SDR exists in \(\Gamma(A, B^{*})\). Let \(\hat{B}\) be a solution of \eqref{P3}. The optimal cost for satisfying each condition in Theorem \ref{t:lin} individually is at most \(\norm{\hat{B}}_w\). Thus, \(\norm{\hat{B}}_w \geq \norm{B^{\prime}}_w\)(obtained from step 1) and \(\norm{\hat{B}}_w \geq \norm{\widetilde{B}}_w\)(obtained from step 2). Since \(\norm{B^{\prime}}_w\geq \norm{B^{\prime\prime}}_w\), we get \(2 \norm{\hat{B}}_w \geq  \norm{B^{\prime\prime}}_w + \norm{\widetilde{B}}_w= \norm{B^{*}}_w\).

Again, note that Algorithm \ref{AlgorithmProblem2} and Algorithm \ref{algorithmstronglyprob2} have polynomial time complexity \(O(dm)\) and \(O((d+m)^3)\) respectively. Rest of the procedure have linear complexity. It is easy to see that Algorithm \ref{algorithmapproximate} has \(O((d+m)^3)\) complexity.
\qed
\\

\bibliographystyle{IEEEtran}
\bibliography{ref}

\bigskip
\bigskip
\end{document}